\documentclass{amsart}

\usepackage{amsfonts}
\usepackage{amssymb}
\usepackage{enumerate}
\usepackage{amsmath}
\usepackage[T1]{fontenc}
\usepackage{graphicx}

\newtheorem{ttt}{Theorem}[section]
\newtheorem{llll}[ttt]{Lemma}
\newtheorem{ccc}[ttt]{Claim}
\newtheorem{eee}[ttt]{Example}
\newtheorem{fff}[ttt]{Fact}
\newtheorem{rrr}[ttt]{Remark}
\newtheorem{sss}[ttt]{Statement}
\newtheorem{ddd}[ttt]{Definition}
\newtheorem{qqq}[ttt]{Question}
\newtheorem{cccc}[ttt]{Corollary}
\newtheorem{nnn}[ttt]{Notation}
\newtheorem{ppp}[ttt]{Problem}
\newtheorem{pppp}[ttt]{Proposition}
\newtheorem{ccccc}[ttt]{Conjecture}

\newcommand{\beq}{\begin{equation} }

\newcommand{\bt}{\begin{ttt}}
\newcommand{\bl}{\begin{llll}}
\newcommand{\bc}{\begin{ccc}}
\newcommand{\bex}{\begin{eee}}
\newcommand{\bfa}{\begin{fff}}
\newcommand{\br}{\begin{rrr}\upshape}
\newcommand{\bst}{\begin{sss}}
\newcommand{\bd}{\begin{ddd}\upshape}
\newcommand{\bq}{\begin{qqq}}
\newcommand{\bnn}{\begin{nnn}}
\newcommand{\bpr}{\begin{ppp}}
\newcommand{\bprop}{\begin{pppp}}
\newcommand{\bcor}{\begin{cccc}}
\newcommand{\bcon}{\begin{ccccc}}

\newcommand{\eeq}{\end{equation}}
\newcommand{\et}{\end{ttt}}
\newcommand{\el}{\end{llll}}
\newcommand{\ec}{\end{ccc}}
\newcommand{\eex}{\end{eee}}
\newcommand{\efa}{\end{fff}}
\newcommand{\er}{\end{rrr}}
\newcommand{\est}{\end{sss}}
\newcommand{\ed}{\end{ddd}}
\newcommand{\eq}{\end{qqq}}
\newcommand{\ecor}{\end{cccc}}
\newcommand{\econ}{\end{ccccc}}
\newcommand{\enn}{\end{nnn}}
\newcommand{\epr}{\end{ppp}}
\newcommand{\eprop}{\end{pppp}}

\newcommand{\bp}{\noindent\textbf{Proof. }}
\newcommand{\ep}{\hspace{\stretch{1}}$\square$\medskip}

\newcommand{\lab}[1]{\label{#1}}

\newcommand{\RR}{\mathbb{R}}

\newcommand{\PP}{\mathbb{P}}

\newcommand{\al}{\alpha}

\newcommand{\de}{\delta}

\newcommand{\om}{\omega}
\newcommand{\si}{\sigma}
\newcommand{\ka}{\kappa}
\newcommand{\la}{\lambda}

\newcommand{\iH}{\mathcal{H}}
\newcommand{\iI}{\mathcal{I}}

\newcommand{\iM}{\mathcal{M}}

\newcommand{\iN}{\mathcal{N}}

\newcommand{\beeq}{\begin{equation}}
\newcommand{\eeeq}{\end{equation}}
\def\su{\subset}

\newcommand{\non}{{\rm non}}
\newcommand{\cof}{{\rm cof}}
\newcommand{\add}{{\rm add}}
\newcommand{\cov}{{\rm cov}}

\newcommand{\diam}{\mathrm{diam}}
\newcommand{\forces}[2]{\Vdash_{#1} \mbox{`` } #2 \mbox{ ''}}

\numberwithin{equation}{section}

\begin{document}

\author[M\'arton Elekes]{M\'arton Elekes$^\ast$}
\address{Alfr\'ed R\'enyi Institute of Mathematics, Hungarian Academy of Sciences,
PO Box 127, 1364 Budapest, Hungary and E\"otv\"os Lor\'and
University, Institute of Mathematics, P\'azm\'any P\'eter s. 1/c,
1117 Budapest, Hungary}
\email{elekes.marton@renyi.mta.hu}
\urladdr{http://www.renyi.hu/$\sim$emarci}
\thanks{$^\ast$Partially supported by the
Hungarian Scientific Foundation grants no.~83726, 104178 and 113047.}

\author[Juris Stepr\=ans]{Juris Stepr\=ans$^\dag$}
\address{Department of Mathematics, York University,
Toronto, Ontario M3J 1P3, Canada}
\email{steprans@mathstat.yorku.ca}
\urladdr{http://www.math.yorku.ca/$\sim$steprans}
\thanks{$^\dag$Supported by a Discovery Grant from NSERC.\\
This research was partially done whilst the authors were visiting fellows at the Isaac Newton Institute for Mathematical Sciences in the programme `Mathematical, Foundational and Computational Aspects of the Higher Infinite' (HIF)}



\title{Set-theoretical problems concerning Hausdorff measures}







\begin{abstract}
J. Zapletal asked if all the forcing notions considered in his monograph
are homogeneous. Specifically, he asked if the forcing consisting of
Borel sets of $\sigma$-finite 2-dimensional Hausdorff measure in
$\mathbb{R}^3$ (ordered under inclusion) is homogeneous. We give a 
partial negative answer to both questions by showing that this $\sigma$-ideal
is not homogeneous.

Let $\mathcal{N}^1_2$ be the $\sigma$-ideal of sets in the plane of 1-dimensional
Hausdorff measure zero. D. H. Fremlin determined the position of the
cardinal invariants of this $\sigma$-ideal in the Cicho\'n Diagram. This required
proving numerous inequalities, and in all but three cases it was known that
the inequalities can be strict in certain models. For one of the remaining
ones Fremlin posed this as an open question in his monograph. We answer
this by showing that consistently $\cov(\mathcal{N}^1_2) >
\cov(\mathcal{N})$, where $\mathcal{N}$ is the usual Lebesgue null
ideal. We also prove that the remaining two inequalities can be strict. 
Moreover, we fit the cardinal invariants
of the $\sigma$-ideal of sets of $\sigma$-finite Hausdorff measure into the diagram.

P. Humke and M. Laczkovich raised the following question. Is it
consistent that there is an ordering of the reals in which all proper
initial segments are Lebesgue null but for every ordering of the reals
there is a proper initial segment that is not null with respect to the
$1/2$-dimensional Hausdorff measure? We determine the values of the
cardinal invariants of the Cicho\'n Diagram as well as the invariants of
the nullsets of Hausdorff measures in the first model mentioned in the
previous paragraph, and as an application we answer this question of
Humke and Laczkovich affirmatively. 
\end{abstract}

\maketitle

\section{Introduction}

Throughout the paper, let $n$ be a positive integer and let $0 < r < n$ be a real number.

\bd
The \emph{$r$-dimensional Hausdorff measure} of a set $H \su \RR^n$ is

\begin{align*}
\mathcal{H}^{r}(H)& = \lim_{\delta\to 0+}\mathcal{H}^{r}_{\delta}(H)
\mbox{, where}\\
\mathcal{H}^{r}_{\delta}(H)& = \inf \left\{ \sum_{k \in \om} (\diam
(A_{k}))^{r} : H \subset \bigcup_{k \in \om} A_{k}, \ \forall k \ \diam (A_k) \le \delta \right\}.
\end{align*}
\ed

\br
\lab{r:capacity}
It is easy to check that $\iH^r(H) = 0$ iff $\iH^r_\infty(H) = 0$.
\er
For more information on this notion see \cite{Fa} or \cite{Mt}.

Let us define the following $\si$-ideal consisting of sets of $\si$-finite $r$-dimensional Hausdorff measure.

\bd
\[
\iI^r_{n, \si-fin} = \{ H \su \RR^n : \exists H_k \su \RR^n, \ \cup_{k \in \om} H_k = H, \ \iH^r(H_k) < \infty \textrm{ for every } k \in \om \}.
\]
\ed

Since it is not hard to see that every set of finite $\iH^r$-measure is contained in a Borel, actually $G_\de$,  set of finite $\iH^r$-measure, this $\sigma$-ideal has a Borel basis (that is, every member of the $\sigma$-ideal is contained in a Borel member of the $\sigma$-ideal).

Following the terminology of \cite{Za} let us define the following notion of forcing.

\bd
\[
\PP_{\iI^r_{n, \si-fin}} = \{B \su \RR^n : B \textrm{ is Borel}, \ B\notin \iI^r_{n, \si-fin}\} \textrm{, ordered under inclusion}.
\]
\ed
For more information on forcing one can also consult \cite{Ku} or \cite{Je}.

In order to be able to formulate our first problem, we need some definitions.

\bd
A notion of forcing $\PP$ is called \emph{homogeneous} if for every $p
\in \PP$ the restriction of $\PP$ below $p$ is forcing equivalent to
$\PP$.
\ed

In his monograph \cite{Za} J. Zapletal poses the following problem.

\bpr 
\lab{q:Z}
\textup{\textbf{(\cite[Question 7.1.3.]{Za})}}
``Prove that some of the forcings presented in this book are not
homogeneous.''
\epr

In fact, we will actually work with the following very closely related notion.

\bd
A  $\sigma$-ideal $\iI$ on a Polish space $X$ is \emph{homogeneous} if for every Borel
set $B \su X$ with $B \notin \iI$ there is a Borel function $f: X \to B$ such that $I \in \iI$ implies $f^{-1}(I) \in \iI$.
\ed

J. Zapletal remarks that ``In all cases encountered in this book the homogeneity of the
forcing and the underlying ideal always come together''. 

Then he also mentions:
``A typical case is that of $\iI$ generated by sets of finite two-dimensional Hausdorff
measure in $\RR^3$.''

In Theorem \ref{t:Z} below we show that this $\sigma$-ideal $\iI^2_{3, \si-fin}$ is indeed non-homogeneous.

\bigskip

Our second problem concerns fitting the cardinal invariants of the $\sigma$-ideal of nullsets of the Hausdorff measures into the Cicho\'n Diagram. For more information on this diagram consult \cite{BJ}.

\bd
Let
\[
\iN^r_n = \{ H \su \RR^n : \iH^r(H) = 0 \}.
\]
\ed

D. H. Fremlin \cite[534B]{Fr} showed that the picture is as follows.

\[
\begin{array}{ccccccccccc}
\cov(\iN) &
\rightarrow &
\cov(\iN^r_n) &
\rightarrow &
\non(\iM) &
\rightarrow &
\cof(\iM)&
\rightarrow &
\cof(\iN^r_n) &
= &
\cof(\iN) \\
 &
 &
 &
 &
\uparrow &
 &
\uparrow &
 &
  \\
\smash{\bigg\uparrow} &
 &
\smash{\bigg\uparrow} &
 &
\mathfrak{b} &
\rightarrow  &
\mathfrak{d} &
 &
\smash{\bigg\uparrow} &
 &
\smash{\bigg\uparrow} \\
 &
 &
 &
 &
\uparrow &
 &
\uparrow &
 &
 \\
\add(\iN) &
= &
\add(\iN^r_n) &
\rightarrow &
\add(\iM)&
\rightarrow &
\cov(\iM)&
\rightarrow &
\non(\iN^r_n) &
\rightarrow &
\non(\iN)
\end{array}
\]

All but three arrows (=inequalities) are known to be strict in
the appropriate models (see e.g. \cite{BJ} for the inequalities not
involving $\iN^r_n$ and \cite{SS} for $\non(\iN^r_n) < \non(\iN)$).
Fremlin, addressing one of these three questions, asked the following.

\bq
\lab{q:F}
\textbf{\textup{(\cite[534Z, Problem (a)]{Fr})}}
Does $\cov(\iN) = \cov(\iN^r_n)$ hold in $ZFC$?
\eq

In Corollary \ref{c:F} below we answer this question in the
negative. The consistent strictness of the remaining two inequaities are
proved in Section \ref{s:further}.

\bigskip

Our last problem was formulated in a recent preprint of P. Humke and
M. Laczkovich \cite{HL}.
Working on certain generalizations of results of Sierpi\'nski and of
Erd\H os they isolated the following definition.

\bd
For a $\sigma$-ideal $\iI$ on $\RR$ let us abbreviate the following statement as 
\[
(*)_{\iI} \iff \textrm{ there exists an ordering of } \RR \textrm{ such that all proper initial segments are in } \iI.
\]
\ed

Using this notation our problem can be formulated as follows.

\bq
\lab{q:HL}
\textbf{\textup{(\cite{HL})}}
Is it consistent that $(*)_{\iN}$ holds but $(*)_{\iN^{1/2}_1}$ fails?
\eq

The following is easy to see and is also shown in \cite{HL}.

\bc
$\add(\iI) = \cov(\iI) \implies (*)_{\iI} \implies \cov(\iI) \le \non(\iI)$.
\ec

Hence it suffices to answer the following question affirmatively.

\bq
\lab{q:HL2}
Is it consistent that $\add(\iN) = \cov(\iN)$ and $\cov(\iN^{1/2}_1) > \non(\iN^{1/2}_1)$?
\eq

In Corollary \ref{c:HL} below we answer this question affirmatively.

\section{Partial answer to Zapletal's question}

\bt
\lab{t:Z}
The $\sigma$-ideal $\iI^2_{3, \si-fin}$ is not homogeneous, partially answering Zapletal's question.
\et

\bp 
Let $B \su \RR^3$ be an arbitrary Borel set with
$\dim_H(B) = \frac52$. Let $f : \RR^3 \to B$ be an arbitrary Borel map.
Then \cite[Theorem 1.4]{Ma} states that for every Borel set $A \su
\RR^n$, Borel map $f: A \to \RR^m$ and $0 \le d \le 1$ there exists a
Borel set $D \su A$ such that $\dim_H D = d \cdot \dim_H A$ and $\dim_H
f(D) \le d \cdot \dim_H f(A)$. Applying this with $n=m=3$, $A=\RR^3$,
and $d = \frac{11}{15}$ we obtain that there exists a Borel set $D \su
\RR^3$ with $\dim_H(D) = \frac{11}5$ such that $\dim_H(f(D)) \le
\frac{11}{15} \cdot \frac52 = \frac{11}6$. Then $\dim_H(D)> 2$ and
$\dim_H(f(D)) < 2$, therefore $f(D) \in \iI^2_{3, \si-fin}$, but $f^{-1}
(f(D)) \supset D \notin \iI^2_{3, \si-fin}$. Since $f$ was arbitrary,
the choice $I = f(D)$ shows that $\iI^2_{3, \si-fin}$ is not
homogeneous.  \ep

\br
The same proof actually yields that for every $0 < r  < n$ the $\sigma$-ideal $\iI^r_{n, \si-fin}$ is not homogeneous.
\er

\section{The model answering the questions of Fremlin and Humke-Laczkovich}

\bl
\lab{l:c=c}
$\cov(\iI^r_{n, \si-fin}) = \cov(\iN^r_n)$
\el

\bp
The inequality $\cov(\iI^r_{n, \si-fin}) \le \cov(\iN^r_n)$ is clear by
$\iN^r_n \su \iI^r_{n, \si-fin}$. In order to prove the opposite
inequality let $\{I_\al\}_{\al < \cov(\iI^r_{n, \si-fin})}$ be a cover
of $\RR^n$ by sets of $\si$-finite $\iH^r$-measure. We can assume that
they are actually of finite $\iH^r$-measure, and also that they are
Borel (even $G_\de$). By the Isomorphism Theorem of Measures
\cite[Thm. 17.41]{Ke} a Borel set of finite $\iH^r$-measure can be
covered by $\cov(\iN)$ many $\iH^r$-nullsets. Therefore $\RR^n$ can be
covered by $\cov(\iI^r_{n, \si-fin}) \cdot \cov(\iN)$ many
$\iH^r$-nullsets. But $\iI^r_{n, \si-fin} \su \iN$ implies
$\cov(\iI^r_{n, \si-fin}) \ge \cov(\iN)$, hence $\RR^n$ can be covered
by $\cov(\iI^r_{n, \si-fin})$ many $\iH^r$-nullsets, proving
$\cov(\iN^r_n) \le \cov(\iI^r_{n, \si-fin})$.
\ep

The following theorem describes the values of all the cardinal invariants of the above diagram in a specific model of $ZFC$.

\bt
\lab{t:main}
It is consistent with $ZFC$ that $\cov(\iN) = \mathfrak{d} = \non(\iN) = \om_1$ and $\cov(\iN^r_n) = \mathfrak{c} = \om_2$.
\et

\bp
Most ingredients of this proof are actually present in \cite{Za}. Let $V$ be a ground model satisfying the Continuum Hypothesis, and let $W$ be obtained by the countable support iteration of $\PP_{\iI^r_{n, \si-fin}}$ of length $\om_2$. Since the forcing $\PP_{\iI^r_{n, \si-fin}}$ is proper by \cite[4.4.2]{Za} and adds a generic real avoiding the Borel members of $\iI^r_{n, \si-fin}$ coded in $V$, we obtain that $\cov(\iI^r_{n, \si-fin}) = \mathfrak{c} = \om_2$ in $W$. Hence, $\cov(\iN^r_n) = \mathfrak{c} = \om_2$ in $W$ by Lemma \ref{l:c=c}.
By \cite[4.4.8]{Za} $\PP_{\iI^r_{n, \si-fin}}$ adds no splitting reals, hence no Random reals, and this is well-known to be preserved by the iteration, thus the Borel nullsets coded in $V$ cover the reals of $W$, therefore $\cov(\iN) = \om_1$ in $W$. Moreover, by \cite[Ex. 3.6.4]{Za} $\iI^r_{n, \si-fin}$ is polar, which is preserved by the iteration, therefore it
preserves outer Lebesgue measure, hence the ground model is not null, thus $\non(\iN) = \om_1$ in $W$.
Finally, the forcing is $\om^\om$-bounding by \cite[4.4.8]{Za}, hence the same holds for the iteration, therefore $\mathfrak{d}=\om_1$ in $W$.
\ep

The following are immediate.

\bcor
\lab{c:F}
Consistently $\cov(\iN) < \cov(\iN^r_n)$, answering Fremlin's question.
\ecor

\bcor
\lab{c:HL}
The answer to Question \ref{q:HL2} is affirmative, hence so is the answer to the question of Humke and Laczkovich.
\ecor

\section{Further results}
\label{s:further}

First, for the sake of completeness, let us now determine the position
of the cardinal invariants of the $\sigma$-ideal $\iI^r_{n, \si-fin}$ in the
diagram.

\bprop
In $ZFC$,
\[
\begin{array}{lll}
\add(\iI^r_{n, \si-fin}) & = & \om_1,\\
\cov(\iI^r_{n, \si-fin}) & = & \cov(\iN^r_n),\\
\non(\iI^r_{n, \si-fin}) & = & \non(\iN^r_n),\\
\cof(\iI^r_{n, \si-fin}) & = & \mathfrak{c}.
\end{array}
\]
\eprop

\bp
Let $\{B_\al\}_{\al < \om_1}$ be a disjoint family of Borel sets of positive finite $\iH^r$-measure, then clearly $\cup_{\al < \om_1} B_\al \notin \iI^r_{n, \si-fin}$ showing $\add(\iI^r_{n, \si-fin}) = \om_1$.

$\cov(\iI^r_{n, \si-fin}) = \cov(\iN^r_n)$ is just Lemma \ref{l:c=c}.

In order to prove $\non(\iI^r_{n, \si-fin}) = \non(\iN^r_n)$, let us assume to the contrary that $\non(\iN^r_n) = \ka < \la = \non(\iI^r_{n, \si-fin})$. Let $H \notin \iN^r_n$ be such that $|H| = \ka$. Then $H$ is of $\si$-finite $\iH^r$-measure, that is $H = \cup_{k \in \om} H_k$ such that $\iH^r(H_k) < \infty $ for every $k \in \om$. Fix $k$ such that $\iH^r(H_k) > 0$. Every set of finite $\iH^r$-measure is contained in a Borel (actually $G_\de$) set of finite $\iH^r$-measure, therefore there exists a Borel set $B \supset H_k$ of positive finite $\iH^r$-measure. Clearly, $|H_k| \le \ka$. By the Isomorphism Theorem of Measures \cite[17.41]{Ke} this implies that $\non(\iN) \le \ka$. But $\iI^r_{n, \si-fin} \su \iN$ yields $\la = \non(\iI^r_{n, \si-fin}) \le \non(\iN) \le \ka$, a contradiction. 

Finally, let $\{B_\al\}_{\al < \mathfrak{c}}$ be a disjoint family of
Borel sets of positive finite $\iH^r$-measure. Since every set of
$\si$-finite $\iH^r$-measure can contain at most countably many of them,
it is easy to see that $\cof(\iI^r_{n, \si-fin}) = \mathfrak{c}$.
\ep

Next we show that the remaining two inequalities in the above extended
Cicho\'n Diagram are also strict in certain models.

Recall that, as usual in set theory, each natural number is identified
with the set of its predecessors, i.e. $k = \{0, \ldots, k-1\}$. Also
recall that $[k]^m = \{A \subset k : |A|=m \}$.

\bt
It is consistent with $ZFC$ that $\cov(\iN^r_n) < \non(\iM)$.
\et

\bp
Let $W$ be the Laver model, that is, the model obtained by iteratively
adding $\omega_2$ Laver reals with countable support over a model $V$
satisfying the Continuum Hypothesis, see \cite{BJ} for the definitions
and basic properties of this model. For example, it is well-known that
$\mathrm{non}(\mathcal{M}) = \omega_2$ in this model.

On the other hand, $W$ satisfies the so called Laver property, an
equivalent form of which is the following:

If $0 < r < n$ and $x\in \prod_{k\in\omega}2^{kn} \cap W$ then there is
\[
T \in \prod_{k\in\omega}[2^{kn}]^{\le 2^{k\frac{r}{2}}} \cap V
\]
such that $x(k)\in T(k)$ for all $k \in \om$. (Here $[A]^{\le d}$ is the system of subsets of $A$ of cardinality at most $d$.) This follows from
\cite[Lemma~6.3.32]{BJ} by letting $f(k) = 2^{kn}$, $S(k) = \{x(k)\}$, and
using and arbitrary positive rational number $s < \frac{r}{2}$.

The following argument takes place in $W$. For every $k \in \om$ let $\psi_k$ be
a bijection from $2^{kn} $ to the set of all cubes of the form
\[
\left[\frac{j_0}{2^k},\frac{j_0+1}{2^k}\right]\times \ldots \times
\left[\frac{j_{n-1}}{2^k},\frac{j_{n-1}+1}{2^k}\right],
\]
where $j_i\in 2^k$ for each $i \in n$. 

For every $T \in \prod_{k\in\omega}[2^{kn}]^{\le 2^{k\frac{r}{2}}}$ define
\[
N_T = \bigcap_{k \in \omega} \bigcup_{j\in T(k)}\psi_k(j).
\]

First we show that $N_T \in \mathcal{N}^r_n$. Note that the diameter of
a cube of side-length $\frac{1}{2^k}$ is
$\sqrt{n}\frac{1}{2^k}$. Clearly, for every $k \in \om$ we have
$\mathcal{H}^r_\infty \left( N_T \right) \le \mathcal{H}^r_\infty \left(
\bigcup_{j\in T(k)}\psi_k(j) \right) \le |T(k)| \left(
\sqrt{n}\frac{1}{2^k}\right)^r \le 2^{k\frac{r}{2}} \left(
\sqrt{n}\frac{1}{2^k}\right)^r  = \sqrt{n}^r 2^{-k\frac{r}{2}}$, which
tends to $0$ as $k$ tends to $\infty$, therefore $\mathcal{H}^r_\infty
\left( N_T \right) = 0$ and consequently, by Remark \ref{r:capacity}, $\mathcal{H}^r \left( N_T
\right) = 0$.

Next we finish the proof by showing that $\{N_T : T \in
\prod_{k\in\omega}[2^{kn}]^{\le 2^{k\frac{r}{2}}} \cap V\}$ is a cover of
$[0,1]^n$ (note that $|V| = \omega_1$ in $W$, and also that if
$\omega_1$ members of $\mathcal{N}^r_n$ cover the unit cube then the
same holds for $\mathbb{R}^n$, hence this implies
$\mathrm{cov}(\mathcal{N}^r_n) = \omega_1$). So let $z\in [0,1]^n$, then
there exists $x\in \prod_{k\in\omega}2^{kn}$ such that $z\in
\psi_k(x(k))$ for each $k \in \om$. Let $T \in
\prod_{k\in\omega}[2^{kn}]^{\le 2^{k\frac{r}{2}}} \cap V$ be such that $x(k)
\in T(k)$ for all $k \in \om$, then it is easy to check that $z \in N_T$,
finishing the proof.
\ep

Next we turn to the consistency of $\cov(\iM) < \non(\iN^r_n)$. First we need some
preparation.

For each $k\in \omega$ let $M_k \in \omega$ be so large that 
\begin{equation}
\label{M_k}
2^k \left( \frac{\sqrt{n}}{M_k} \right)^r < \frac{1}{2^k}. 
\end{equation}

\bd
Let $C_k$ be the set of all cubes of the form 
\[
\left[\frac{j_0}{M_k},\frac{j_0+1}{M_k}\right]\times \ldots \times \left[\frac{j_{n-1}}{M_k},\frac{j_{n-1}+1}{M_k}\right],
\]
where $j_i \in M_k$ for each $i \in n$. 
Let $\mathbb C_k$ consist of all sets that can be written as the union of $2^k$ elements of $C_k$.
\ed

\bl
For every partition $\mathbb C_k = \bigcup_{i \in 2^k} X_i$ there is some $i \in 2^k$ such that $\cup X_i = [0,1]^n$.
\el

\bp
Otherwise, pick $x_i \notin \cup X_i$ and cubes $Q_i \in C_k$ containing $x_i$, then $\bigcup_{i \in 2^k} Q_i \in \mathbb{C}_k$ belongs to one of the $X_i$, yielding a contradiction.
\ep

\bd
Now we define the \emph{norm function} $\nu \colon \bigcup_{k \in \omega} \mathcal{P}(\mathbb{C}_k) \to \omega$ as follows.
For $X\subset \mathbb C_k$ define $\nu(X) \geq 1$ if $\cup X = [0,1]^n$ and
define $\nu(X) \geq  j+1$ if for every partition $X= X_0\cup X_1$ there
is $i\in 2$ such that $\nu(X_i)\geq j$.
\ed

\bl
$\nu(\mathbb C_k) \geq k+1$.
\el

\begin{proof}
Otherwise, we could iteratively split $\mathbb C_k$ into pieces so that at stage $m$ we have a partition into $2^m$ many sets each with norm at most $k-m$, hence eventually we could have a partition into $2^k$ many sets none of which covers $[0,1]^n$, contradicting the previous lemma. 
\end{proof}

\bl
\label{l:thin}
If $X\subset \mathbb C_k$ and $\nu(X) \geq j$ and $y\in [0,1]^n$ then
$\nu( \{ H\in X : y\in H \} ) \geq j-1$.
\el

\begin{proof}
We may assume $j>1$. Let $X_0 = \{ H\in X : y\in H \}$ and $X_1 =
 \{ H \in X : y \notin H \}$. Then either $\nu(X_0)\geq j-1\geq 1$ or
 $\nu(X_1)\geq j-1\geq 1$. But note that $\nu(X_1) \not\geq  1$ since $y\notin \cup X_1$.
\end{proof}

In this paper a \emph{finite
sequence} will mean a function defined on a natural number, the \emph{length}
of the sequence $t$, denoted by $|t|$ is simply $\mathrm{dom}(t)$.
Moreover, a \emph{tree} will mean a set of finite sequences closed under
initial segments. Then for $t,s \in T$ we have
$t \subset s$ iff $s$ end-extends $t$ and this partial order is indeed a
tree in the usual sense. For a $t \in T$ let us denote by $\mathrm{succ}_T (t)$ the
set of immediate successors of $t$ in $T$.

\medskip
Now let us define the following forcing notion. 

\bd
Let $T \in \mathbb{P}$ iff
\begin{enumerate}[(1)]
\item
$T$ is a non-empty tree,
\item
\label{2}
for every $t \in T$ and $k < |t|$ we have $t(k) \in \mathbb{C}_k$,
\item
for every $t \in T$ we have $\mathrm{succ}_T(t) \neq \emptyset$,
\item
for every $t \in T$ there exists $s \in T$, $s \supset t$ with
$|\mathrm{succ}_T(s)| > 1$,
\item
for every $K \in \omega$ the set $\{t \in T : |\mathrm{succ}_T(t)| > 1
\textrm{ and } \nu(\mathrm{succ}_T(t)) \le K\}$ is finite.
\end{enumerate}

If $T, T' \in \mathbb{P}$ then define
\[
T \le_\mathbb{P} T' \iff T \subset T'.
\]
We will usually simply write $\le$ for $\le_\mathbb{P}$. Clearly,
$1_\mathbb{P}$ is the set of all finite sequences satisfying (\ref{2}).
\ed

\br
A $t \in T$ with $|\mathrm{succ}_T(t)| > 1$ is called a \emph{branching node}.
For $t \in T$ define $T[t] = \{s \in T : s \subset t \textrm{ or } s \supset
t\}$. It is easy to see that if $t \in T \in \mathbb{P}$ then $T[t] \in
\mathbb{P}$ and $T[t] \le T$.
\er

\br
Forcing notions of this type are discussed in
\cite{RS} in great generality. However, in order to keep the paper relatively self-contained we
also include the rather standard proofs here, but note that most of the
techniques below can already be found in \cite{RS}.
\er

\bl
\label{l:proper}
$\mathbb{P}$ is proper.
\el

\begin{proof}
Let $\mathfrak{M}$ be a countable elementary submodel, and recall that $T \in
\mathbb{P}$ is \emph{$\mathfrak{M}$-generic} if for every dense open subset $D
\subset \mathbb{P}$ with $D \in \mathfrak{M}$ we have $T \forces{}{\dot{G}
\cap D \cap \mathfrak{M} \neq \emptyset}$, where $\dot{G}$ is a name for the generic filter.
Also recall that properness means that whenever a condition $T \in
\mathfrak{M}$ is given then there exists an $\mathfrak{M}$-generic $T' \le
T$. We construct this $T'$ by a so called fusion argument. 

Let the sequence $D_0, D_1, \ldots$ enumerate the dense open subsets of
$\mathbb{P}$ that are in $\mathfrak{M}$. During the construction we make sure
that all objects we pick ($t$, $s$, $t_s'$, $s'$, $L_k$, $L_k^+$ $S_s$, etc.) are in $\mathfrak{M}$. The whole construction, and
hence $T'$, will typically not be in $\mathfrak{M}$.

We define the set of branching notes
of $T'$ `level-by level' as follows. Let $t \in T$ be a branching
node with $\nu(\mathrm{succ}_T(t)) > 0$ and set $L_0 =\{t\}$. Also define $L_0^+ = \mathrm{succ}_T(t)$. Moreover,
for every $s \in L_0^+$ also fix a $S_s \le T[s]$ with $S_s \in D_0$ (this is
possible, since $D_0$ is dense). This finishes the $0$th step of the fusion.

Now, if $L_k, L_k^+$, and for every $s \in L_k^+$ a condition $S_s \le T[s]$
have already been defined then for every $s \in L_k^+$ we pick a $t_s' \in S_s$
with $\nu(\mathrm{succ}_{S_s}(t_s')) > k+1$. Let $L_{k+1} = \{t_s' : s \in
L_k^+\}$, and define $L_{k+1}^+ = \bigcup_{s \in L_k^+}
\mathrm{succ}_{S_s}(t_s')$. Now, for every $s' \in L_{k+1}^+$ pick a $S_{s'}
\le S_s[s']$ with $S_{s'} \in D_{k+1}$. This finishes the
$k+1$st step of the fusion.

Finally, define $T'$ as the closure of $\bigcup_{k \in \omega} L_k$ under
initial segments (this is the same as the closure of $\bigcup_{k \in \omega} L_k^+$ under
initial segments). It is easy to check that $T' \in \mathbb{P}$ and $T' \le
T$. It remains to show that $T'$ is $\mathfrak{M}$-generic. So let $k \in
\omega$ be fixed, and we need to show that $T' \forces{}{\dot{G}
\cap D_k \cap \mathfrak{M} \neq \emptyset}$.

Before the proof let us make three remarks. First, it is easy to see from the
construction that if $T'' \le T'$ then for every $k \in \omega$ there exists
$s \in L_k^+ \cap T''$. Second, it can also be seen from the construction that
$T'[s] \le S_s$ for every $s$. Third, if $S \in D \cap \mathfrak{M}$ then obviously $S \forces{}{\dot{G}
\cap D \cap \mathfrak{M} \neq \emptyset}$, since $S \forces{}{S \in \dot{G}
\cap D \cap \mathfrak{M}}$.

Now we prove $T' \forces{}{\dot{G} \cap D_k \cap \mathfrak{M} \neq
  \emptyset}$. We prove this by showing that for every $T'' \le T'$ there
exists $T''' \le T''$ forcing this. Let $T'' \le T'$ be given. Then, by the
above remark there exists $s \in L_k^+ \cap T''$. Set $T''' = T''[s]$, then
clearly $T''' \le T''$. Finally, $T''' = T''[s] \le T'[s] \le S_s \in D_k \cap
\mathfrak{M}$, hence $S_s \forces{}{\dot{G}
\cap D_k \cap \mathfrak{M} \neq \emptyset}$, hence $T'''$ forces the same, finishing the proof.
\end{proof}

\bl
\label{l:limsup}
If $(H_k)_{k \in \omega} \in \Pi_{k \in \omega} \mathbb{C}_k$ then
$\mathcal{H}^r( \bigcap_{m \in \omega} \bigcup_{k \ge m} H_k) = 0$.
\el

\begin{proof}
For every $m \in \om$, using \eqref{M_k}, we have 
\[
\mathcal{H}^r_\infty \left( \bigcap_{m \in \omega} \bigcup_{k \ge m} H_k \right)
\le
\mathcal{H}^r_\infty \left( \bigcup_{k \ge m} H_k \right)
\le
\sum_{k \ge m} \mathcal{H}^r_\infty \left( H_k \right)
\le
\]
\[
\le
\sum_{k \ge m} 2^k \left( \frac{\sqrt{n}}{M_k}\right)^r
\le
\sum_{k \ge m} \frac{1}{2^k}
=
\frac{1}{2^{m-1}},
\]
hence $\mathcal{H}^r_\infty ( \bigcap_{m \in \omega} \bigcup_{k \ge m} H_k) =
0$, therefore $\mathcal{H}^r( \bigcap_{m \in \omega} \bigcup_{k \ge m} H_k ) = 0$.
\end{proof}

\br
In the usual way, by slight abuse of notation, the generic filter $G$ can be
thought of as a sequence $G = (G_k)_{k \in \omega} \in \Pi_{k \in \omega}
\mathbb{C}_k$. What we will formally need is that if a generic filter $G$ is given, then $\bigcap_{T \in G} T$ defines such
a sequence, hence $G_k$ makes sense.
\er

\bl
\label{l:main}
If $G$ is a generic filter over a ground model $V$ then $V[G] \models V \cap
[0,1]^n \subset \bigcap_{m \in \omega} \bigcup_{k \ge m} G_k$.
\el

\begin{proof}
Fix $y \in V \cap [0,1]^n$. In order to show that $1_\mathbb{P} \forces{}{V \cap
[0,1]^n \subset \bigcap_{m \in \omega} \bigcup_{k \ge m} \dot{G}_k}$ we show that for every $T \in \mathbb{P}$ there is $T' \le T$ forcing this. So let $T$ be given, and define $T'$ as follows. Starting from the root of $T$, we recursively thin out $T$ such that for every $t \in T$ with $\nu(\mathrm{succ}_T (t)) \ge 1$ we cut off all the nodes $s \in \mathrm{succ}_T (t)$ with $y \notin s(|s|-1)$. One can easily check using Lemma \ref{l:thin} that $T' \in \mathbb{P}$ and $T' \le T$. So it suffices to show that for every $m \in \omega$ we have $T' \forces{}{y \in \bigcup_{k \ge m} \dot{G}_k}$. Hence let $T'' \le T'$ be given, we need to find $T''' \le T''$ forcing this. Pick $t \in T''$ with $|t| \ge m$ and $\nu(\mathrm{succ}_{T''} (t)) \ge 1$. This implies that the successors of $t$ were thinned out, hence $y \in s(|s|-1)$ for every $s \in \mathrm{succ}_{T''} (t)$. Fix such an $s$, and define $T''' = T''[s]$. Then $T''' = T''[s] \forces{}{\dot{G}_{|s|-1} = s(|s|-1) \ni y}$, finishing the proof.
\end{proof}

\bl
$\PP$ is $\omega^\omega$-bounding.
\el

\begin{proof}
For $f, g \in \omega^\omega$ we write $f \le g$ if $f(n) \le g(n)$ for every $n \in \omega$. Let $\dot{f} \in \omega^\omega$ be a name. We claim that $1_\mathbb{P} \forces{}{\exists g \in V \cap \omega^\omega \textrm{ such that } \dot{f} \le g}$. It suffices to show that for every $T$ there exists $T' \le T$ and $g \in V \cap \omega^\omega$ such that $T' \forces{}{\dot{f} \le g}$. We will construct this $T'$ by a fusion argument similar to that of Lemma \ref{l:proper}. 

Let $t \in T$ be a node with $\nu(\mathrm{succ}_T(t)) > 0$ and set $L_0 =\{t\}$.
Also define $L_0^+ = \mathrm{succ}_T(t)$. Moreover,
for every $s \in L_0^+$ also fix a $S_s \le T[s]$ and $m_s \in \omega$ such that $S_s \forces{}{\dot{f} (0) = m_s}$ (this is possible by the basic properties of forcing).
This finishes the $0$th step of the fusion.

Now, if $L_k, L_k^+$, and for every $s \in L_k^+$ a condition $S_s \le T[s]$
have already been defined then for every $s \in L_k^+$ we pick a $t_s' \in S_s$
with $\nu(\mathrm{succ}_{S_s}(t_s')) > k+1$. Let $L_{k+1} = \{t_s' : s \in
L_k^+\}$, and define $L_{k+1}^+ = \bigcup_{s \in L_k^+}
\mathrm{succ}_{S_s}(t_s')$. Now, for every $s' \in L_{k+1}^+$ pick a $S_{s'}
\le S_s[s']$ and $m_{s'} \in \omega$ with $S_{s'} \forces{}{\dot{f} (k+1) = m_{s'}}$. This finishes the
$k+1$st step of the fusion.

Finally, define $T'$ as the closure of $\bigcup_{k \in \omega} L_k^+$ under
initial segments. It is easy to check that $T' \in \mathbb{P}$ and $T' \le
T$. Define $g(k) = \max\{ m_s : s \in L_k^+\}$ (the maximum exists, since this set is finite). 
It remains to show that $T' \forces{}{\dot{f} \le g}$. So let $k \in
\omega$ be fixed, and let $T'' \le T'$ be given. Pick $s \in T'' \cap L_k^+$, and define $T''' = T''[s]$. Then $T''' = T''[s] \le S_s \forces{}{ \dot{f} (k) = m_s }$, hence $T''' \forces{}{ \dot{f} (k) \le g(k)}$, finishing the proof.
\end{proof}

\bt
It is consistent with $ZFC$ that $\cov(\iM) < \non(\iN^r_n)$.
\et

\begin{proof}
Let $V$ be a model satisfying the Continuum Hypothesis, and let
$V_{\omega_2}$ be the model obtained by an $\omega_2$-long countable
support iteration of $\PP$. Let $(V_\alpha)_{\alpha \le \omega_2}$
denote the intermediate models. Since $\PP$ is proper and adds a real, by
standard arguments the continuum is $\om_2$ in $V_{\om_2}$.

On the one hand, $\mathbb{P}$ is $\omega^\omega$-bounding, hence so is
its iteration. Therefore the iteration adds no Cohen reals, hence the
meagre Borel sets coded in $V$ cover $V_{\omega_2} \cap \mathbb{R}^n$,
hence $\mathrm{cov} (\mathcal{M}) = \omega_1$.

On the other hand, if $H \in V_{\omega_2}, |H| = \omega_1$ then, by a
standard reflection argument, $H \subset V_\alpha$ for some $\alpha <
\omega_2$. Hence, by Lemma \ref{l:main} and Lemma \ref{l:limsup} we have
$V_{\alpha+1} \models \mathcal{H}^r (H) = 0$. Therefore, since
$\mathcal{H}^r (H) = 0$ means the existence of certain covers, and by
absoluteness the corresponding covers exist in $V_{\omega_2}$, we obtain
$V_{\omega_2} \models \mathcal{H}^r (H) = 0$. Hence, $\mathrm{non}
(\mathcal{N}^r_n) = \omega_2$. Therefore the proof
is complete.
\end{proof}

\section{Open problems}

First we reiterate Zapletal's question \cite[Question 7.1.3.]{Za} in its original form.

\bpr
Is the forcing notion $\PP_{\iI^2_{3, \sigma-fin}}$ homogeneous?
\epr

Moreover, Theorem \ref{t:Z} also leaves open the following.

\bpr
Does there exist a Borel set $B \subset \RR^3$ with $B \notin \iI^2_{3, \sigma-fin}$ such that $\iI^2_{3, \sigma-fin}$ is not homogeneous below $B$?
\epr

\bpr
Let $0 < r < s < n$. Does $\non(\iN^r_n) = \non(\iN^s_n)$ hold in $ZFC$?
\epr

Analogously,

\bpr
Let $0 < r < s < n$. Does $\cov(\iN^r_n) = \cov(\iN^s_n)$ hold in $ZFC$?
\epr

\section*{Acknowledgement}

The authors would like to thank P. Frenkel for spotting a disturbing error in an earlier version of the paper.

\end{document}